\input epsf
\documentclass{amsart}
\usepackage{amssymb,latexsym,hyperref}

\newcommand{\doublesubscript}[3]{
\displaystyle\mathop{\displaystyle #1_{#2}}_{#3} }

\numberwithin{equation}{section}
\newtheorem{theorem}{Theorem}[section]
\newtheorem{proposition}[theorem]{Proposition}

\newtheorem{definition}[theorem]{Definition}

\newtheorem{remark}[theorem]{Remark}
\newtheorem{lemma}[theorem]{Lemma}

\def\CC{\mathbb{C}}
\def\FF{\mathbb{F}}
\def\gg{\mathfrak{g}}
\def\hh{\mathfrak{h}}
\def\ii{\mathbf{i}}
\def\l{\ell}
\def\nn{\mathfrak{n}}
\def\RR{\mathbb{R}}
\def\veps{\varepsilon}
\def\ZZ{\mathbb{Z}}

\newcommand{\mat}[4]{\left(\!\!\begin{array}{cc}
#1 & #2 \\ #3 & #4 \\
\end{array}\!\!\right)}

\begin{document}

\title{Connected components of real double Bruhat cells}

\author{Andrei Zelevinsky}
\address{\noindent Department of Mathematics, Northeastern University,
  Boston, MA 02115}
\email{andrei@neu.edu}

\thanks{The author's research was supported in part
by NSF grant \#DMS-9971362.}

\date{April 24, 2000}

\maketitle

\section{Introduction}
\label{sec:introduction}

The main geometric objects of study in this paper are \emph{double
Bruhat cells} $G^{u,v} = B u B  \cap B_- v B_-$ in a
simply-connected connected complex semisimple group $G$; here $B$ and $B_-$
are two opposite Borel subgroups in $G$, and $u$ and $v$ any two
elements of the Weyl group $W$. Double Bruhat cells were
introduced and studied in \cite{FZ} as a geometric framework for
the study of total positivity in semisimple groups; they are also
closely related to symplectic leaves in the corresponding
Poisson-Lie groups (see \cite{DCP,HKKR}). It will be convenient for us
to replace $G^{u,v}$ with a reduced double Bruhat cell $L^{u,v}$
introduced in \cite{BZ99}. The variety $L^{u,v}$ can be identified
with the quotient of $G^{u,v}$ modulo the left (or right) action
of the maximal torus $H = B \cap B_-$.

As shown in \cite{FZ,BZ99}, an algebraic variety $L^{u,v}$ is
biregularly isomorphic to a Zariski open subset of an affine space
of dimension $m = \l(u)+\l(v)$, where $\l(u)$ is the length of $u$
in the Coxeter group $W$. However, the smooth topology of $L^{u,v}$
can be quite complicated. A first step towards
understanding this topology is enumerating the connected
components of the real part $L^{u,v}(\RR)$. In the
case when $G$ is simply-laced, a conjectural answer was given in
\cite[Conjecture~4.1]{SSVZ}. Here we prove this conjecture
and extend the result to an arbitrary semisimple group $G$.
The answer is given in the following terms: as shown in \cite{SSVZ}
for $G$ simply-laced, every reduced word $\ii$ of $(u,v) \in W \times W$
gives rise to a subgroup $\Gamma_\ii (\FF_2) \subset GL_m (\FF_2)$
generated by symplectic transvections (here $\FF_2$ is the
$2$-element field). We extend the construction of $\Gamma_\ii
(\FF_2)$ to an arbitrary $G$ (it is still generated by
transvections but not necessarily by symplectic ones). Extending
\cite[Conjecture~4.1]{SSVZ}, we show that the connected components of
$L^{u,v}(\RR)$ are in a natural bijection with the
$\Gamma_\ii (\FF_2)$-orbits in $\FF_2^{m}$.
As explained in \cite{SSVZ}, this provides a
far-reaching generalization of the results in \cite{SSV1,SSV2};
this also refines and generalizes results in \cite{RI,RI2}.

Our proof uses methods and results developed in \cite{FZ,BZ99}.
First, it was shown there that every reduced word $\ii$ of $(u,v)
\in W \times W$ gives rise to a biregular isomorphism between the
complex torus $\CC_{\neq 0}^m$ and a Zariski open subset $U_\ii
\subset L^{u,v}$. We refine this result by showing that the
complement $L^{u,v} - U_\ii$ is the union of $m$ divisors
$\{M_{k,\ii} = 0\}$, where $M_{1,\ii}, \dots, M_{m,\ii}$ are some
irreducible regular functions on $L^{u,v}$. We further show that
every $\ii$-bounded index $n \in [1,m]$ (see
Section~\ref{sec:background} for the definition) gives rise to a
regular function $M'_{n,\ii}$ on $L^{u,v}$ such that replacing the
divisor $\{M_{n,\ii} = 0\}$ with $\{M'_{n,\ii} = 0\}$ leads to
another ``toric chart" $U_{n,\ii}$ in $L^{u,v}$. Then we prove
that the connected components of the real part of the union of
charts $U_\ii \  \bigcup \ \bigcup_n U_{n,\ii}$ are in a natural
bijection with the $\Gamma_\ii (\FF_2)$-orbits in $\FF_2^{m}$.
Finally, we show that the complement in $L^{u,v}$ of this union of
charts has codimension $\geq 2$, so the connected components of
$L^{u,v}(\RR)$ are enumerated in the same way as those of the real
part of $U_\ii \ \bigcup \ \bigcup_n U_{n,\ii}$.

According to \cite{FZ,BZ99}, each $M_{k,\ii}$ is a ``twisted
(generalized) minor" on $G$. We show that each $M'_{n,\ii}$ is
obtained by the same twist from a regular function on $G$ which is
no longer a minor but can be expressed as a sum of two Laurent
monomials in minors. These regular functions are of
independent interest for the study of the dual canonical bases in the
ring of regular functions $\CC[G]$ and its $q$-deformation.

The paper is organized as follows.
After recalling the necessary background, we formulate
our main result (Theorem~\ref{th:components}) in
Section~\ref{sec:background}.
In Section~\ref{sec:lemma}, we formulate a lemma (Lemma~\ref{lem:main})
that plays the crucial role
in our proof of Theorem~\ref{th:components}, and then show how
this lemma implies the theorem.
The proof of Lemma~\ref{lem:main} is given in
Section~\ref{sec:proof lemma}.
Finally, Section~\ref{sec:examples} discusses some examples and
applications of the results in Sections~\ref{sec:lemma} and
\ref{sec:proof lemma}.

\section{Main theorem}
\label{sec:background}

To formulate our main result, let us recall the necessary
background from \cite{FZ,SSVZ,BZ99}. Let $G$ be a simply connected
semisimple algebraic group with the Dynkin graph $\Pi$. Let $B$
and $B_-$ be two $\RR$-split opposite Borel subgroups, $N$ and
$N_-$ their unipotent radicals, $H = B \cap B_-$ an $\RR$-split
maximal torus of $G$, and $W = {\rm Norm}_G (H) / H$ the Weyl
group of $G$. Let $\gg = {\rm Lie}(G)$ be the Lie algebra of $G$,
and $\hh = {\rm Lie}(H)$ the Cartan subalgebra of $\gg$. Let
$\{\alpha_i: i \in \Pi\}$ be the system of simple roots in $\hh^*$
for which the corresponding root subgroups are contained in~$N$.
Let $\{\alpha_i^\vee: i \in \Pi\}$ be the corresponding system of
simple coroots in $\hh$, and $A = (a_{ij} = \alpha_j
(\alpha_i^\vee))$ be the \emph{Cartan matrix}. Thus, for $i \neq
j$ the indices $i$ and $j$ are adjacent in $\Pi$ if and only if
$a_{ij} a_{ji} \neq 0$; we shall denote this by $\{i,j\} \in \Pi$.
For every $i \in \Pi$, let $\varphi_i: SL_2 \to G$ denote the
corresponding canonical $SL_2$-embedding.

The Weyl group $W$ is canonically identified with the Coxeter
group $W(A)$ generated by the involutions $s_i$ for $i \in \Pi$
subject to the relations $(s_i s_j)^{d_{ij}} = e$ for all $i \neq
j$, where $d_{ij} =2$ (resp.\ $3,4,6$) if $a_{ij}a_{ji} = 0$
(resp.\ $1,2,3$).
A word $\ii = (i_1, \ldots, i_m)$ in the alphabet $\Pi$ is a
\emph{reduced word} for $w \in W$ if $w = s_{i_1} \cdots s_{i_m}$,
and $m$ is the smallest length of such a factorization. The length
$m$ of any reduced word for $w$ is called the \emph{length} of $w$
and denoted by $m = \l (w)$. Let $R(w)$ denote the set of all
reduced words for $w$. The identification $W = W(A)$ is given by
$s_i = \overline {s_i} H$, where $$\overline {s_i} = \varphi_i
\mat{0}{-1}{1}{0} \in {\rm Norm}_G (H) \ .$$ The representatives
$\overline {s_i} \in G$ satisfy the braid relations $\overline
{s_i} \overline {s_j} \overline {s_i} \cdots = \overline {s_j}
\overline {s_i} \overline {s_j} \cdots$ (with $d_{ij}$ factors on
each side); thus, the representative $\overline w$ can be
unambiguously defined for any $w \in W$ by requiring that
$\overline {uv} = \overline {u} \cdot \overline {v}$ whenever $\l
(uv) = \l (u) + \l (v)$.

The ``double" group $W \times W$ is also a Coxeter group. The
corresponding graph $\tilde \Pi$ is the union of two disconnected
copies of $\Pi$. We identify the vertex set of $\tilde \Pi$ with
$\{+1, -1\} \times \Pi$, and write a vertex $(\pm 1, i) \in \tilde
\Pi$ simply as $\pm i$. For each $i \in \Pi$, we set
$\veps (\pm i) = \pm 1$ and $|\pm i| = i$. Thus, two
vertices $i$ and $j$ of $\tilde \Pi$ are adjacent if and only if
$\veps (i) = \veps (j)$ and $\{|i|, |j|\} \in \Pi$. In this
notation, a reduced word for a pair $(u,v) \in W \times W$ is an
arbitrary shuffle of a reduced word for $u$ written in the
alphabet $-\Pi$ and a reduced word for $v$ written in the alphabet
$\Pi$.

The group $G$ has two \emph{Bruhat decompositions}, with respect
to $B$ and $B_-\,$: $$G = \bigcup_{u \in W} B u B = \bigcup_{v \in
W} B_- v B_-  \ . $$ The \emph{double Bruhat cells}~$G^{u,v}$ are
defined by $G^{u,v} = B u B  \cap B_- v B_- \,$.

Following \cite{BZ99}, we define the \emph{reduced double Bruhat
cell} $L^{u,v} \subset G^{u,v}$ as follows:
\begin{equation}
\label{eq: left sections} L^{u,v} = N \overline u N  \cap B_- v
B_- \ .
\end{equation}
The maximal torus $H$ acts freely on $G^{u,v}$ by left (or right)
translations, and $L^{u,v}$ is a section of this action. Thus,
$G^{u,v}$ is biregularly isomorphic to $H \times L^{u,v}$, and all
properties of $G^{u,v}$ can be translated in a straightforward way
into the corresponding properties of $L^{u,v}$ (and vice versa).
In particular, Theorem~1.1 in \cite{FZ} implies that $L^{u,v}$ is
biregularly isomorphic to a Zariski open subset of an affine space
of dimension $\l(u)+\l(v)$.

The \emph{real part} of $G$ is the subgroup $G(\RR)$ of $G$ generated by all the
subgroups $\varphi_i (SL_2(\RR))$. For any subset $L \subset G$, we
define its real part by $L(\RR) = L \cap G(\RR)$.

Now let us fix a pair $(u,v) \in W \times W$, and
let $m = \l (u) + \l(v)$. Let $\ii = (i_1, \ldots, i_m) \in
R(u,v)$ be any reduced word for $(u,v)$.
We associate to $\ii$ an $m \times m$ matrix $(C_{kl})$
in the following way: set $C_{kl} = 1$ if $|i_k| = |i_l|$
and $C_{kl} = - a_{|i_k|,|i_l|}$ if $|i_k| \neq |i_l|$.

Following \cite{SSVZ}, we
associate with $\ii$ a directed graph $\Sigma (\ii)$ on the set of
vertices $[1,m] = \{1, 2, \ldots, m\}$. For $l \in [1,m]$, we
denote by $l^- = l^-_\ii$ the maximal index $k$ such that $1 \leq
k < l$ and $|i_k| = |i_l|$; if $|i_k| \neq |i_l|$ for $1 \leq k <
l$ then we set $l^- = 0$. The edges of $\Sigma(\ii)$ are now
defined as follows.

\begin{definition}
\label{def:edges} {\rm A pair $\{k,l\} \subset [1,m]$ with $k<l$
is an edge of $\Sigma(\ii)$ if it satisfies one of the following
three conditions:

(i) $k=l^-$;

(ii) $k^- < l^- < k$, $\{|i_k|, |i_l|\} \in \Pi$, and
$\veps(i_{l^-})=\veps(i_{k})$;

(iii) $l^-< k^- < k$, $\{|i_k|, |i_l|\} \in \Pi$, and
$\veps(i_{k^-})=-\veps(i_{k})$.

\noindent The edges of type (i) are called \emph{horizontal}, and
those of types (ii) and (iii) \emph{inclined}. A horizontal (resp.
inclined) edge $\{k,l\}$ with $k < l$ is directed from $k$ to $l$
if and only if $\veps (i_{k}) = +1$ (resp. $\veps (i_{k}) = -1$).
We shall write $(k \to l) \in \Sigma(\ii)$ if $k \to l$ is a
directed edge of $\Sigma(\ii)$.}
\end{definition}

We now associate to each $n \in [1,m]$ a transvection $\tau_n =
\tau_{n, \ii}: \ZZ^m \to \ZZ^m$ defined as follows: if $\tau_n
(\xi_1, \dots, \xi_m) = (\xi'_1, \dots, \xi'_m)$ then $\xi'_k =
\xi_k$ for $k \neq n$, and
\begin{equation}
\label{eq:transvection}
\xi'_n  = \xi_n  - \sum_{(k \to n) \in
\Sigma(\ii)}C_{kn} \xi_k + \sum_{(n \to l) \in \Sigma(\ii)} C_{ln} \xi_l
\end{equation}
(note that if $G$ is simply-laced then all the coefficients
$C_{kn}$ and $C_{ln}$ in (\ref{eq:transvection})are equal to $1$,
so (\ref{eq:transvection}) becomes formula (2.4) in \cite{SSVZ}).
We call an index $n \in
[1,m]$ \emph{$\ii$-bounded} if $n^- > 0$. Let $\Gamma_{\ii}$
denote the group of linear transformations of $\ZZ^m$ generated by
the transvections $\tau_n$ for all $\ii$-bounded indices $n \in
[1,m]$. Let $\Gamma_{\ii}(\FF_2)$ denote the group of linear
transformations of the $\FF_2$-vector space $\FF_2^m$ obtained
from $\Gamma_{\ii}$ by reduction modulo $2$ (recall that $\FF_2$
is the $2$-element field).

We are finally ready to formulate our main result.

\begin{theorem}
\label{th:components} For every reduced word $\ii \in R(u,v)$,
the connected components of $L^{u,v}(\RR)$ are in a natural bijection
with the $\Gamma_\ii (\FF_2)$-orbits in $\FF_2^{m}$.
\end{theorem}

Note that in Theorem~\ref{th:components} we only need the modulo
$2$ reductions of transvections $\tau_n$, so the formula
(\ref{eq:transvection}) could be simplified as follows:
$$\xi'_n
= \xi_n  + \sum_{(k,n) \in \Sigma(\ii)}C_{kn} \xi_k \ .$$
We prefer the form (\ref{eq:transvection}) because it is suggested by
the construction of toric charts in $L^{u,v}$ which is our main
ingredient in proving Theorem~\ref{th:components}.

\section{Main lemma}
\label{sec:lemma}

As before, let $G$ be a simply
connected connected complex semisimple group with the Dynkin graph $\Pi$.
We fix a pair $(u,v) \in W \times W$, let $m = \l (u) + \l(v)$,
and fix a reduced word $\ii = (i_1, \ldots, i_m) \in R(u,v)$.

\begin{lemma}
\label{lem:main} There exist regular functions $M_1, \dots, M_m$
on $L^{u,v}$ with the following properties:

\smallskip

\noindent (1) If $k \in [1,m]$ is not $\ii$-bounded then $M_k$
vanishes nowhere on $L^{u,v}$.

\smallskip

\noindent (2) The map $(M_1, \dots, M_m): L^{u,v} \to
\CC^m$ restricts to a biregular isomorphism $U_\ii \to \CC_{\neq
0}^m$, where $U_\ii$ is the locus of all $x \in L^{u,v}$ such that
$M_k(x) \neq 0$ for all $k \in [1,m]$.

\smallskip

\noindent (3) For every $\ii$-bounded $n \in [1,m]$, the rational
function $M'_k$ defined by
\begin{equation}
\label{eq:M'} M'_{n} M_{n} = \prod_{(k \to n) \in \Sigma_\ii}
M_{k}^{C_{kn}} + \prod_{(n \to l) \in \Sigma_\ii} M_{l}^{C_{ln}}
\end{equation}
is regular on $L^{u,v}$.

\smallskip

\noindent (4) For every $\ii$-bounded $n \in [1,m]$,
the map $(M_1, \dots, M_{n-1}, M'_n, M_{n+1}, \dots, M_m): L^{u,v} \to \CC^m$
restricts to a biregular isomorphism $U_{n,\ii} \to \CC_{\neq 0}^m$, where
$U_{n,\ii}$ is the locus of all $x \in L^{u,v}$ such that $M'_n (x)
\neq 0$ and $M_k (x) \neq 0$ for all $k \in [1,m] - \{n\}$.

\smallskip

\noindent (5) The functions $M_k$ and $M'_n$ take real values on
$L^{u,v} (\RR)$, and the biregular isomorphisms in (2) and (4)
restrict to biregular isomorphisms $U_{\ii}(\RR) \to \RR_{\neq
0}^m$ and $U_{n,\ii}(\RR) \to \RR_{\neq 0}^m$.
\end{lemma}

The functions $M_k = M_{k,\ii}$ in Lemma~\ref{lem:main} were
introduced in \cite[(4.13)]{BZ99}.
We recall the definition and prove Lemma~\ref{lem:main}
in the next section; in the
rest of this section we show that it implies
Theorem~\ref{th:components}. To be more precise, we shall prove
that the bijection in Theorem~\ref{th:components} can be defined
as follows. For every $\xi = (\xi_1, \dots, \xi_m) \in \FF_2^{m}$,
let $U_\ii (\xi)$ denote the set of all $x \in U_\ii (\RR)$ such
that $(-1)^{\xi_k} M_k (x) > 0$ for all $k$. For every $Y \subset
L^{u,v} (\RR)$, let $\overline Y$ denote the closure of $Y$ in
$L^{u,v} (\RR)$ in the real topology.

\begin{theorem}
\label{th:components explicit} The correspondence $\Omega \mapsto
\bigcup_{\xi \in \Omega} \overline {U_\ii (\xi)}$ is a bijection
between $\Gamma_\ii (\FF_2)$-orbits in $\FF_2^{m}$ and connected
components of $L^{u,v}(\RR)$.
\end{theorem}

We split the proof of Theorem~\ref{th:components explicit} into
several lemmas. Let us abbreviate $X = L^{u,v}$, and let $\CC[X]$
be the ring of regular functions on $X$. Since $X$ is isomorphic
to a Zariski open subset of $\CC^m$, the ring $\CC[X]$ is a unique
factorization domain. By property (1) in Lemma~\ref{lem:main}, if
$k$ is not $\ii$-bounded then $M_k$ is an invertible element of
$\CC[X]$.

\begin{lemma}
\label{lem:regular monomials} A Laurent monomial $P = M_1^{d_1}
\cdots M_m^{d_m}$ is a regular function on $X$ if and only if $d_n
\geq 0$ for any $\ii$-bounded $n$.
\end{lemma}

\proof The ``if" part is trivial. To prove the ``only if" part,
fix an $\ii$-bounded index $n$, and consider the restriction of
$P$ to the Zariski open subset $U_{n,\ii} \subset X$. By property
(4) in Lemma~\ref{lem:main}, if $P \in \CC[X]$ then $M_n^{d_n}$ is
a regular function on $U_{n,\ii}$ and so it must be a Laurent
polynomial in $M_1, \dots, M_{n-1}, M'_n, M_{n+1}, \dots, M_m$. In
view of (\ref{eq:M'}), this implies that $d_n \geq 0$, as desired.
\endproof

\begin{lemma}
\label{lem:Mk irreducible} For every $\ii$-bounded $n$, the
function $M_n$ is an irreducible element of $\CC[X]$.
\end{lemma}

\proof Notice that every $P \in \CC[X]$ restricts to a regular
function on the Zariski open subset $U_{\ii} \subset X$. By
property (2) in Lemma~\ref{lem:main}, $P$ is a Laurent polynomial
in $M_1, \dots, M_m$. It follows that if $M_n$ is the product of
two regular functions $P$ and $Q$ then both $P$ and $Q$ must be
Laurent monomials in $M_1, \dots, M_m$. By Lemma~\ref{lem:regular
monomials}, one of the factors $P$ and $Q$ must be  a Laurent
monomial in the variables $M_k$ for $k$ not $\ii$-bounded, hence
is an invertible element of $\CC[X]$. Therefore, $M_n$ is
irreducible.
\endproof

\begin{lemma}
\label{lem:M'k irreducible} For every $\ii$-bounded $n$, the
function $M'_n$ is equal to some irreducible element $M''_n \in
\CC[X]$ times a Laurent monomial in $M_1, \dots, M_{n-1}, M_{n+1},
\dots, M_m$.
\end{lemma}

\proof Let $P \in \CC[X]$ be an irreducible factor of $M'_n$.
Restricting $P$ to $U_{n,\ii}$, we conclude that $P$ is a Laurent
monomial in $M_1, \dots, M_{n-1}, M'_n, M_{n+1}, \dots, M_m$.
Restricting $P$ to $U_{\ii}$ and using property (2) in
Lemma~\ref{lem:main}, we see that $P$ must be also a Laurent
polynomial in $M_1, \dots, M_m$. By (\ref{eq:M'}), this implies
that the exponent of $M'_n$ in $P$ written as a Laurent monomial in $M_1, \dots,
M_{n-1}, M'_n, M_{n+1}, \dots, M_m$ must be nonnegative.
It follows that there is an irreducible
factor $M''_n$ of $M'_n$ which is equal to $M'_n$ times a Laurent
monomial in $M_1, \dots, M_{n-1}, M_{n+1}, \dots, M_m$, while the
rest of the factors are just Laurent monomials in $M_1, \dots,
M_{n-1}, M_{n+1}, \dots, M_m$.
\endproof

We set $U = U_\ii \ \bigcup \ \bigcup_n U_{n,\ii}$.

\begin{lemma}
\label{lem:complement} The complement $X - U$ is the locus of all
$x \in X$ such that $M_n (x) = M_k (x)= 0$ for two distinct
$\ii$-bounded indices $n$ and $k$, or $M_n (x) = M''_n (x) = 0$
for some $\ii$-bounded $n$. The variety $X - U$ has (complex)
codimension $\geq 2$ in $X$.
\end{lemma}

\proof Suppose $x \in X - U$. Since $x \notin U_\ii$, property (1)
in Lemma~\ref{lem:main} implies that $M_n (x) = 0$ for some
$\ii$-bounded $n$. Since $x \notin U_{n,\ii}$, it follows that
either $M''_n (x) = 0$, or $M_k (x) = 0$ for some $\ii$-bounded
$k \neq n$.

The converse inclusion is obvious. Finally, the statement that $X
- U$ has codimension $\geq 2$ in $X$ is clear since $X - U$ is the
union of finitely many subvarieties, each given by two (distinct)
irreducible equations.
\endproof

Now consider the real part $U(\RR)  = U_\ii (\RR) \ \bigcup \
\bigcup_n U_{n,\ii}(\RR)$. By Lemma~\ref{lem:complement} and
property (5) in Lemma~\ref{lem:main}, the complement $X(\RR) -
U(\RR)$ has real codimension $\geq 2$ in $X(\RR)$. Therefore, the
connected components of $X(\RR)$ (in the real topology) are
closures of the connected components of $U(\RR)$. It remains to
show that Theorem~\ref{th:components explicit} holds with $X(\RR)$
replaced by $U(\RR)$. For a subset $Y \subset U(\RR)$ we now
denote by $\overline Y$ the closure of $Y$ in $U(\RR)$. The
role of transvections $\tau_n$ is explained by the following
lemma.

\begin{lemma}
\label{lem:neighbors} Let $\xi^{(1)}$ and $\xi^{(2)}$ be two
distinct vectors in $\FF_2^m$. Then $\overline {U(\xi^{(1)})} \cap
\overline {U(\xi^{(2)})} \neq \emptyset$ if and only if $\xi^{(2)}
= \tau_n (\xi^{(1)})$ for some $\ii$-bounded index $n$.
\end{lemma}

\proof Suppose $x \in U(\RR)$ belongs to the intersection
$\overline {U(\xi^{(1)})} \cap \overline {U(\xi^{(2)})}$. Then
$M_k (x) = 0$ whenever $\xi^{(1)}_k \neq \xi^{(2)}_k$. Using
Lemma~\ref{lem:complement}, we see that there is a unique $n$ such
that $\xi^{(1)}_n \neq \xi^{(2)}_n$; furthermore, this index $n$
is $\ii$-bounded, and $M'_n (x) \neq 0$. Since any neighborhood of
$x$ intersects both $U(\xi^{(1)})$ and $U(\xi^{(2)})$, it follows
that the two monomials on the right hand side of (\ref{eq:M'})
must have opposite signs at $x$. Let us write $\xi_k = \xi^{(1)}_k
= \xi^{(2)}_k$ for $k \neq n$. Then we have
$$\xi^{(2)}_n -
\xi^{(1)}_n = 1 = \sum_{(n \to l) \in \Sigma(\ii)} C_{ln} \xi_l -
\sum_{(k \to n) \in \Sigma(\ii)}C_{kn} \xi_k \ .$$
Comparing this with (\ref{eq:transvection}), we conclude that
$\xi^{(2)} = \tau_n (\xi^{(1)})$, as claimed.

Conversely, suppose $\xi^{(2)} = \tau_n (\xi^{(1)}) \neq
\xi^{(1)}$, and let $\xi_k = \xi^{(1)}_k = \xi^{(2)}_k$ for $k \neq n$.
Then $$\sum_{(n \to l) \in \Sigma(\ii)} C_{ln} \xi_l \neq
\sum_{(k \to n) \in \Sigma(\ii)}C_{kn} \xi_k \ .$$
This implies
that there exists a point $x \in U_{n,\ii}(\RR)$ such that
$(-1)^{\xi_k} M_k (x) > 0$ for all $k \neq n$, and the right hand
side of (\ref{eq:M'}) vanishes at $x$. Then any neighborhood of
$x$ contains points with the signs of all $M_k$ for $k \neq n$
unchanged and with the right hand side of (\ref{eq:M'} positive
(as well as negative). Thus, $x \in \overline {U(\xi^{(1)})} \cap
\overline {U(\xi^{(2)})}$, and we are done.
\endproof

Now we are ready to complete the proof of
Theorem~\ref{th:components explicit}. Let $\Omega$ be a
$\Gamma_\ii (\FF_2)$-orbit in $\FF_2^{m}$, and consider the
corresponding closed subset $Y_\Omega = \bigcup_{\xi \in \Omega}
\overline {U_\ii(\xi)}$ of $U(\RR)$. Each $U_\ii(\xi)$ is
a copy of $\RR_{> 0}^m$ and so is connected. Using the ``if" part of
Lemma~\ref{lem:neighbors}, we conclude that $Y_\Omega$ is
connected (since the closure of a
connected set and the union of two non-disjoint
connected sets are connected as well).
On the other hand, by the ``only if" part of the same lemma, all
the sets $Y_\Omega$ are pairwise disjoint. Thus, they are the
connected components of $U(\RR)$, and we are done.
\endproof

\section{Proof of Lemma~\ref{lem:main}}
\label{sec:proof lemma}

\subsection{The functions $M_k$}
\label{sec:Mk}
We start by recalling the definition of the functions $M_k = M_{k,\ii}$
given in \cite[(4.13)]{BZ99}.
First of all, recall that the \emph{weight lattice} $P$ of $G$
can be thought of as the group of rational multiplicative characters of~$H$
written in the exponential notation:
a weight $\gamma\in P$ acts by $a \mapsto a^\gamma$.
The lattice $P$ is also identified with the additive group of all
$\gamma \in \hh^*$ such that
$\gamma (\alpha_i^\vee) \in \ZZ$ for all $i \in \Pi$.
Thus, $P$ has a $\ZZ$-basis $\{\omega_i: i \in \Pi\}$  of \emph{fundamental weights}
given by $\omega_j (\alpha_i^\vee) = \delta_{i,j}$.

We now recall from \cite{FZ} the definition of \emph{generalized minors}.
Denote by $G_0=N_- H N$ the open subset of elements $x\in G$ that
have Gaussian decomposition; this (unique) decomposition will be written as
$x = [x]_- [x]_0 [x]_+ \,$.
For $u,v \in W$ and $i \in \Pi$, the \emph{(generalized) minor}
$\Delta_{u \omega_i, v \omega_i}$
is the regular function on $G$ whose restriction to the open set
${\overline {u}} G_0 {\overline {v}}^{-1}$ is given by
\begin{equation}
\label{eq:Delta-general}
\Delta_{u \omega_i, v \omega_i} (x) =
(\left[{\overline {u}}^{\ -1}
   x \overline v\right]_0)^{\omega_i} \ .
\end{equation}
As shown in \cite{FZ}, $\Delta_{u \omega_i, v \omega_i}$ depends on
the weights $u \omega_i$ and $v \omega_i$ alone, not on the particular
choice of $u$ and~$v$.
It is easy to see that the generalized minors are distinct
irreducible elements of the ring $\CC[G]$ of regular functions on $G$.
In the special case $G=SL_n\,$, the generalized minors are nothing but
the ordinary minors of a matrix.

According to \cite[Proposition~4.3]{BZ99},
an element $x \in G^{u,v}$ belongs to $L^{u,v}$ if and only if
$[{\overline u}^{\ -1} x]_0 = 1$, or equivalently if
$\Delta_{u \omega_i, \omega_i} (x) = 1$ for any $i \in  [1,r]$.

We fix a pair $(u,v) \in W \times W$ and a double reduced word
$\ii = (i_1, \ldots, i_m)\in R(u,v)$.
Recall from \cite[Definition~4.6 and Theorem~4.7]{BZ99} that there is
a biregular isomorphism $\psi^{u,v}$ between $L^{u,v}$ and
$L^{v, u}$ given by
\begin{equation}
\label{eq:psi-u,v-x}
\psi^{u,v}(x) =
[(\overline{v} x^{\iota})^{-1}]_+ \, \overline{v} \,
([\overline{u}^{\ -1} x]_+)^{\iota} \ ;
\end{equation}
here $x \mapsto x^\iota$ is the involutive antiautomorphism of $G$
given by
\begin{equation}
\label{eq:iota}
\varphi_i \mat{a}{b}{c}{d}^\iota = \varphi_i \mat{d}{b}{c}{a} \ .
\end{equation}

Recall that the length $m$ of $\ii$ is equal to $\l (u) + \l (v)$.
For $k\in [1,m]$, denote
\begin{equation}
\label{eq:v_leq k}
u_{\geq k} = \doublesubscript{\prod}{l = m, \ldots, k}
{i_l \in - \Pi}
s_{|i_l|} \ , \quad
v_{<k} = \doublesubscript{\prod}{l = 1, \ldots, k-1}
{i_l \in \Pi } s_{i_l} \ .
\end{equation}
This notation means that in the first (resp.\ second) product in
(\ref{eq:v_leq k}), the index $l$ is decreasing (resp. increasing);
for example, if $\Pi = \{1,2,3\}$ and
$\ii= (-2, 1, -3, 3, 2, -1, -2, 1, -1)$,
then, say, $u_{\geq 7} = s_1 s_2 \,$ and $v_{<7} = s_1 s_3 s_2\,$.

Following  \cite[(4.13)]{BZ99}, we define a regular function $M_k = M_{k,\ii}$ on $L^{u,v}$ by
\begin{equation}
\label{eq:M-factors}
M_k (x) = M_{k,\ii} (x) =
\Delta_{v_{<k} \omega_{|i_k|}, u_{\geq k} \omega_{|i_k|}} (\psi^{u,v} (x)) \ .
\end{equation}

\subsection{Properties (1), (2) and (5)}
\label{sec:1-2}
To prove property (1) in Lemma~\ref{lem:main},
notice that if $k$ is not $\ii$-bounded and $|i_k| = i$ then
(\ref{eq:M-factors}) turns into
$M_{k,\ii} (x) = \Delta_{\omega_{i}, u^{-1} \omega_{i}} (\psi^{u,v}(x))$.
It remains to show that $\Delta_{\omega_{i}, u^{-1} \omega_{i}}$
vanishes nowhere on $L^{v,u}$.
In fact, a stronger statement holds: $\Delta_{\omega_{i}, u^{-1} \omega_{i}}$
vanishes nowhere on $B_- u B_-$.
This follows from the definition (\ref{eq:Delta-general}) and the well-known
inclusion $B_- u B_- u^{-1} \subset G_0$ (cf.~\cite[Proposition 2.10]{FZ}).

As for property (2) in Lemma~\ref{lem:main}, it follows from the
solution to the so-called factorization problem given in
\cite[Theorem 1.9]{FZ} (or rather from its modification in
\cite[Theorem 4.8]{BZ99}).
To formulate it, we need some notation.

For every $i \in \Pi$ and $t \in \CC_{\neq 0}$, we denote
$$x_i (t) = \varphi_i \mat{1}{t}{0}{1}, \,\, y_i (t) = \varphi_i \mat{1}{0}{t}{1},
\, \, t^{\alpha_i^\vee} = \varphi_i \mat{t}{0}{0}{t^{-1}} \ ;$$
following \cite{BZ99}, we also denote
\begin{equation}
\label{eq:xnegative}
x_{- i} (t) = y_i (t) t^{- \alpha_i^\vee}  = \varphi_i \mat{t^{-1}}{0}{1}{t} \, .
\end{equation}
For any word $\ii= (i_1, \ldots, i_m)$ in the alphabet $\tilde \Pi$,
let us define the \emph{product map} $x_\ii: \CC_{\neq 0}^m \to G$ by
\begin{equation}
\label{eq:productmap}
x_\ii (t_1, \ldots, t_m) =  x_{i_1} (t_1) \cdots x_{i_m} (t_m) \ .
\end{equation}
For $k\in [1,m]$, we denote $k^+ = \min\{l:l>k, |i_l| = |i_k|\}$,
so that $k^+$ is the next occurrence of an index
$\pm i_k$ in $\ii$; if $k$ is the last occurrence of $\pm i_k$ in $\ii$ then we set
$k^+ = m+1$.
We also adopt the convention that $M_{m+1} (x) = 1$.

The following reformulation of Theorem~4.8 in \cite{BZ99} provides
a refinement of property (2) in Lemma~\ref{lem:main}.

\begin{theorem}
\label{th:t-through-x}
Let $\ii = (i_1, \ldots, i_m)$ be a double reduced word for $(u,v)$, and let
$(M_1, \dots, M_m) \in \CC_{\neq 0}^m$.
Then there is a unique $x \in L^{u,v}$ such that $M_{k,\ii} (x) = M_k$ for
$k \in [1,m]$.
This element $x$ has the form $x= x_\ii (t_1, \ldots, t_m)$,
with the factorization parameters $t_k$ given by:
if $i_k \in - \Pi$ then
\begin{equation}
\label{eq:negative t-through-x}
t_k =  M_k/M_{k^+} \ ;
\end{equation}
if $i_k \in \Pi$ then
\begin{equation}
\label{eq:positive t-through-x}
t_k =  \frac{1}{M_k M_{k^+}} \prod_{l: \ l^- < k < l} M_l^{- a_{|i_l|, i_k}} \ .
\end{equation}
\end{theorem}

\begin{remark}
{\rm We see that the parameters $t_1,\dots,t_m$ in the factorization $x= x_\ii (t_1, \ldots, t_m)$
are related to $M_1, \dots, M_m$ by an invertible monomial transformation.
The inverse of this monomial transformation can be computed explicitly: a direct calculation shows that
\begin{equation}
\label{eq:inverse-monomial}
M_k = M_{k,\ii} (x) = \prod_{l \geq k} t_l^{- \veps (i_l) v_{< l}^{-1} v_{<k} \omega_{|i_k|}
(\alpha^\vee_{|i_l|})} \ .
\end{equation}
}
\end{remark}

Finally, property (5) in Lemma~\ref{lem:main}
is clear since each $M_k$ is just a Laurent monomial in the
factorization parameters $t_1, \dots, t_m$, while each $M'_n$ is
the sum of two Laurent monomials; therefore they take real values
when all $t_k$ are real.

\subsection{Property (3)}
\label{sec:regularity}
To prove property (3) in Lemma~\ref{lem:main}, we shall construct
a new family of regular functions on the whole group $G$.
Let $\ii = (i_1, \dots, i_m)$ be a reduced word for $(u,v) \in W
\times W$ such that $|i_1| = |i_m| = i$ for some $i \in \Pi$, and
$|i_k| \neq i$ for $1 < k < m$.
Let $E_{\pm} = \{k \in [2,m-1] \ : \ \veps (i_k) = \pm 1 \}$, and
let $J_\pm = \{i\} \cup \{|i_k| \ : \ k \in E_\pm, k^+ = m+1\} \subset \Pi$.
Let $\Delta' = \Delta'_\ii$ be the rational function on $G$ defined
by one of the following four equations.

\smallskip

\noindent {\bf Case 1.} If $i_1 = i_m = i$ then
\begin{eqnarray*}
\Delta' \Delta_{s_i \omega_i, \omega_i} = \Delta_{\omega_i, \omega_i}
 \doublesubscript {\prod}{k \in E_+}{k^+ \in E_- \cup \{m+1\}}
 \Delta_{v_{\leq k} \omega_{i_k}, u_{> k} \omega_{i_k}}^{- a_{i_k,i}}
 \\[.1in]
 + \Delta_{v \omega_i, \omega_i}
 \doublesubscript {\prod}{k \in E_+}{k^- \in E_- \cup \{0\}}
 \Delta_{v_{< k} \omega_{i_k}, u_{> k} \omega_{i_k}}^{- a_{i_k,i}} \ .
 \end{eqnarray*}

\smallskip

\noindent {\bf Case 2.} If $i_1 = i$ and $i_m = -i$ then
\begin{eqnarray*}
\Delta' \Delta_{s_i \omega_i, s_i \omega_i} = \Delta_{\omega_i, s_i \omega_i}
\Delta_{s_i \omega_i, \omega_i}
 \doublesubscript {\prod}{k \in E_+}{k^+ \in E_-}
 \Delta_{v_{\leq k} \omega_{i_k}, u_{> k} \omega_{i_k}}^{- a_{i_k,i}}
\\[.1in]
+  \doublesubscript {\prod}{k \in E_+}{k^- \in E_- \cup \{0\}}
 \Delta_{v_{< k} \omega_{i_k}, u_{> k} \omega_{i_k}}^{- a_{i_k,i}}
 \prod_{j \in \Pi - J_+} \Delta_{v \omega_j, \omega_j}^{- a_{ji}}
 \ .
 \end{eqnarray*}

 \smallskip

\noindent {\bf Case 3.} If $i_1 = -i$ and $i_m = i$ then
 \begin{eqnarray*}
\Delta' \Delta_{\omega_i, \omega_i} = \Delta_{\omega_i, u^{-1} \omega_i}
\Delta_{v \omega_i, \omega_i}
 \doublesubscript {\prod}{k \in E_+}{k^- \in E_-}
 \Delta_{v_{< k} \omega_{i_k}, u_{> k} \omega_{i_k}}^{- a_{i_k,i}}
\\[.1in]
+  \doublesubscript {\prod}{k \in E_-}{k^- \in E_+ \cup \{0\}}
 \Delta_{v_{< k} \omega_{|i_k|}, u_{\geq k} \omega_{|i_k|}}^{- a_{|i_k|,i}}
 \prod_{j \in \Pi - J_-} \Delta_{v \omega_j, \omega_j}^{- a_{ji}}
 \ .
 \end{eqnarray*}

 \smallskip

\noindent {\bf Case 4.} If $i_1 = i_m = -i$ then
\begin{eqnarray*}
\Delta' \Delta_{\omega_i, s_i \omega_i} = \Delta_{\omega_i, u^{-1} \omega_i}
 \doublesubscript {\prod}{k \in E_-}{k^+ \in E_+ \cup \{m+1\}}
 \Delta_{v_{< k} \omega_{|i_k|}, u_{> k} \omega_{|i_k|}}^{- a_{|i_k|,i}}
 \\[.1in]
 + \Delta_{\omega_i, \omega_i}
 \doublesubscript {\prod}{k \in E_-}{k^- \in E_+ \cup \{0\}}
 \Delta_{v_{< k} \omega_{|i_k|}, u_{\geq k} \omega_{|i_k|}}^{- a_{|i_k|,i}} \ .
 \end{eqnarray*}

 \smallskip

\begin{theorem}
\label{th:delta' regular}
In each of the above four cases, $\Delta' = \Delta'_\ii$ is a
regular function on $G$.
\end{theorem}

Before proving Theorem~\ref{th:delta' regular}, we show that it implies
property (3) in Lemma~\ref{lem:main}.
Let $(u,v)$ be an arbitrary pair of elements of $W$,
and fix a reduced word $\ii = (i_1, \ldots, i_m) \in R(u,v)$.
Let $n$ be an $\ii$-bounded index in $[1,m]$, and let $\ii'$
denote the subword $(i_{n^-}, \dots, i_n)$ of $\ii$.
We claim that the rational function $M'_n$ on $L^{u,v}$ defined by
(\ref{eq:M'}) is given by
\begin{equation}
\label{eq:M' through delta'}
M'_n (x) = M'_{n,\ii} (x) =
\Delta'_{\ii'}(\overline {u_{> n}}^{\ -1} \psi^{u,v} (x) \overline {v_{< n^-}}) \ ,
\end{equation}
and so is regular.
To see this, let us evaluate the defining equation
for $\Delta'_{\ii'}$ at the point
$\overline {u_{> n}}^{\ -1} \psi^{u,v} (x) \overline {v_{< n^-}}$.
Remembering the definition (\ref{eq:Delta-general}) of generalized
minors, and the definition (\ref{eq:M-factors}) of the functions
$M_k$, a direct check shows that, in each of the above four cases, the
corresponding equality turns into the equation (\ref{eq:M'}) with
$M'_n$ given by (\ref{eq:M' through delta'}).

It remains to prove Theorem~\ref{th:delta' regular}.
Our main tool will be the following identity established in
\cite[Theorem~1.17]{FZ}:
\begin{equation}
\label{eq:minors-Dodgson}
\Delta_{v' \omega_i, u' \omega_i} \Delta_{v' s_i \omega_i, u' s_i \omega_i}
- \Delta_{v's_i \omega_i, u' \omega_i} \Delta_{v' \omega_i, u' s_i \omega_i}
= \prod_{j \in \Pi - \{i\}} \Delta_{u' \omega_j, v' \omega_j}^{- a_{ji}}
\end{equation}
for any $u' ,v' \in W$ and $i \in \Pi$
such that $\l (u's_i) = \l (u') + 1$ and $\l (v's_i) = \l (v') + 1$.

To prove that $\Delta'_\ii$ is regular on $G$, we first consider
the case when $\ii$ is ``non-mixed," i.e., $k < l$ for each
$k \in E_-$ and $l \in E_+$.
Then the defining equation for $\Delta' = \Delta'_\ii$ simplifies as
follows.
Denote $S_{\pm} = \{|i_k|: k \in E_{\pm}\} \subset \Pi$.

\smallskip

\noindent {\bf Case 1 (non-mixed).} If $i_1 = i_m = i$ then
\begin{eqnarray*}
\Delta' \Delta_{s_i \omega_i, \omega_i} = \Delta_{\omega_i, \omega_i}
\prod_{j \in S_+} \Delta_{v \omega_j, \omega_j}^{- a_{ji}} + \Delta_{v \omega_i, \omega_i}
\prod_{j \in S_+} \Delta_{\omega_j,  \omega_j}^{- a_{ji}} \ .
 \end{eqnarray*}

\smallskip

\noindent {\bf Case 2 (non-mixed).} If $i_1 = i$ and $i_m = -i$ then
\begin{eqnarray*}
\Delta' \Delta_{s_i \omega_i, s_i \omega_i} = \Delta_{\omega_i, s_i \omega_i}
\Delta_{s_i \omega_i, \omega_i} +  \prod_{j \in \Pi - \{i\}} \Delta_{\omega_j, \omega_j}^{- a_{ji}}
 \ .
 \end{eqnarray*}

 \smallskip

\noindent {\bf Case 3 (non-mixed).} If $i_1 = -i$ and $i_m = i$ then
 \begin{eqnarray*}
\Delta' \Delta_{\omega_i, \omega_i} = \Delta_{\omega_i, u^{-1} \omega_i}
\Delta_{v \omega_i, \omega_i} \prod_{j \in S_+ \cap S_-}
 \Delta_{\omega_j, \omega_j}^{- a_{ji}}  \\[.1in]
+  \prod_{j \in S_-} \Delta_{\omega_j, u^{-1} \omega_j}^{- a_{ji}}
 \prod_{j \in \Pi - (S_- - S_+)} \Delta_{v \omega_j, \omega_j}^{- a_{ji}}
 \ .
 \end{eqnarray*}

 \smallskip

 \noindent {\bf Case 4 (non-mixed).} If $i_1 = i_m = -i$ then
\begin{eqnarray*}
\Delta' \Delta_{\omega_i, s_i \omega_i} = \Delta_{\omega_i, u^{-1} \omega_i}
\prod_{j \in S_-} \Delta_{\omega_j, u^{-1} \omega_{j}}^{- a_{ji}}
+ \Delta_{\omega_i, \omega_i} \prod_{j \in S_-}
 \Delta_{\omega_{j}, u^{-1} \omega_{j}}^{- a_{ji}} \ .
 \end{eqnarray*}

 \smallskip

By (\ref{eq:minors-Dodgson}), in Case 2 we have
$\Delta' = \Delta_{\omega_i, \omega_i}$.
Cases 1 and 4 are equivalent to each other in view of the identity
$\Delta_{\gamma, \delta} (x^T) = \Delta_{\delta, \gamma}(x)$,
where $x \mapsto x^T$ is the involutive antiautomorphism of $G$
given by
\begin{equation*}
\varphi_i \mat{a}{b}{c}{d}^T = \varphi_i \mat{a}{c}{b}{d}
\end{equation*}
(see \cite[Proposition~2.7]{FZ}).
It remains to show that $\Delta'$ is regular in each of the cases 1 and 3.

Let us start with Case 3.
Multiplying both monomials on the right hand side of the
corresponding equation with the monomial
$$\prod_{j \in \Pi -  \{i\} - (S_+ \cap S_-)} \Delta_{\omega_j, \omega_j}^{- a_{ji}}
= \prod_{j \in \Pi - \{i\} - S_-} \Delta_{\omega_j, u^{-1} \omega_j}^{- a_{ji}}
 \prod_{j \in S_- - S_+} \Delta_{v \omega_j, \omega_j}^{- a_{ji}}$$
and using (\ref{eq:minors-Dodgson}), we obtain
$$\Delta_{\omega_i, u^{-1} \omega_i} \Delta_{v \omega_i, \omega_i}
\prod_{j \in \Pi - \{i\}} \Delta_{\omega_j, \omega_j}^{- a_{ji}}+
\prod_{j \in \Pi - \{i\}} \Delta_{\omega_j, u^{-1} \omega_j}^{- a_{ji}}
 \prod_{j \in \Pi - \{i\}} \Delta_{v \omega_j, \omega_j}^{- a_{ji}}$$
$$= \Delta_{\omega_i, u^{-1} \omega_i} \Delta_{v \omega_i, \omega_i}
(\Delta_{\omega_i, \omega_i} \Delta_{s_i \omega_i, s_i \omega_i}
- \Delta_{s_i \omega_i, \omega_i} \Delta_{\omega_i, s_i \omega_i})$$
$$+ (\Delta_{\omega_i, \omega_i} \Delta_{s_i \omega_i, u^{-1} \omega_i}
- \Delta_{s_i \omega_i, \omega_i} \Delta_{\omega_i, u^{-1} \omega_i})
(\Delta_{\omega_i, \omega_i} \Delta_{v \omega_i, s_i \omega_i}
- \Delta_{v \omega_i, \omega_i} \Delta_{\omega_i, s_i \omega_i})$$
$$= \Delta_{\omega_i, \omega_i} \det
\left(\begin{array}{ccc}
\Delta_{\omega_i, u^{-1} \omega_i} & \Delta_{\omega_i, s_i \omega_i} &\Delta_{\omega_i, \omega_i}\\
\Delta_{s_i\omega_i, u^{-1} \omega_i} & \Delta_{s_i\omega_i, s_i \omega_i} &\Delta_{s_i\omega_i, \omega_i}\\
0 & \Delta_{v \omega_i, s_i \omega_i} &\Delta_{v \omega_i, \omega_i}\\
\end{array}\right) \ .$$
Since all ``principal minors" $\Delta_{\omega_j, \omega_j}$ are distinct
irreducible elements of $\CC[G]$, it follows that $\Delta_{\omega_i, \omega_i}$ is
relatively prime with
$\prod_{j \in \Pi -  \{i\} - (S_+ \cap S_-)} \Delta_{\omega_j, \omega_j}^{-a_{ji}}$.
Therefore,
$$\Delta' = \det
\left(\begin{array}{ccc}
\Delta_{\omega_i, u^{-1} \omega_i} & \Delta_{\omega_i, s_i \omega_i} &\Delta_{\omega_i, \omega_i}\\
\Delta_{s_i\omega_i, u^{-1} \omega_i} & \Delta_{s_i\omega_i, s_i \omega_i} &\Delta_{s_i\omega_i, \omega_i}\\
0 & \Delta_{v \omega_i, s_i \omega_i} &\Delta_{v \omega_i, \omega_i}\\
\end{array}\right) /
\prod_{j \in \Pi -  \{i\} - (S_+ \cap S_-)} \Delta_{\omega_j,
\omega_j}^{-a_{ji}}$$
is a regular function on $G$, as required.

The argument in Case 1 is similar (and simpler).
Let us only give the final answer: the function $\Delta'$ is now
given by
$$\Delta' = (\Delta_{\omega_i, \omega_i} \Delta_{v \omega_i, s_i \omega_i}
- \Delta_{v \omega_i, \omega_i} \Delta_{\omega_i, s_i \omega_i}) /
\prod_{j \in \Pi -  \{i\} - S_+} \Delta_{\omega_j, \omega_j}^{-a_{ji}} \ ,$$
and it is again a regular function on $G$, as required.

We shall deduce the general case in Theorem~\ref{th:delta' regular} from
the non-mixed case just considered.
Note that every reduced word $\ii$ in each of the cases  1 -- 4 is
obtained from the corresponding non-mixed word by a sequence of
$2$-\emph{moves} each of which interchanges a pair of consecutive
indices $i_k$ and $i_{k+1}$ with $k \in E_-$ and $k + 1 \in E_+$.
It suffices to show that if $\ii'$ is obtained from $\ii$ by such
a move then the regularity of $\Delta'_\ii$ implies that of
$\Delta'_{\ii'}$.
We shall only treat Case 1; the argument in the other three cases is the same.

Let $P_1$ and $P_2$ (resp. $P'_1$ and $P'_2$) be two monomials in
the right hand side of the defining equation for $\Delta'_{\ii}$
(resp. for $\Delta'_{\ii'}$).
Thus, we assume that $P_1 + P_2$ is divisible by $\Delta_{s_i \omega_i, \omega_i}$
in $\CC[G]$, and need to show that the same is true for $P'_1 + P'_2$.
Let us abbreviate $v' = v_{< k}$ and $u' = u_{> k+1}$, where $k$
and $k+1$ are two positions involved in the $2$-move that turns
$\ii$ into $\ii'$.
It is clear from the definitions that $P'_1 = P_1$ and $P'_2 = P_2$
unless $- i_k = i_{k+1} = j$ for some $j \in \Pi$ such that $a_{ji} < 0$.
In the latter case, we have
$$P'_1 = P_1 \left(\frac{\Delta_{v' s_j \omega_j, u' s_j \omega_j}}
{\Delta_{v'\omega_j, u' s_j \omega_j}^{\delta_1}
\Delta_{v' s_j \omega_j, u'\omega_j}^{\delta_2}}\right)^{-a_{ji}}, \
P'_2 = P_2 \left(\frac{\Delta_{v'\omega_j, u' s_j \omega_j}^{1- \delta_1}
\Delta_{v' s_j \omega_j, u'\omega_j}^{1- \delta_2}}
{\Delta_{v'\omega_j, u'\omega_j}}\right)^{-a_{ji}},$$
where $\delta_1 = 1$ (resp. $\delta_2 = 0$) if $k^- > 1$ and
$i_{k^-} = j$ (resp. $(k+1)^+ < m$ and $i_{(k+1)^+} = j$),
otherwise $\delta_1 = 0$ (resp. $\delta_2 = 1$).
Since the common denominator
$(\Delta_{v'\omega_j, u' s_j \omega_j}^{\delta_1}
\Delta_{v' s_j \omega_j, u'\omega_j}^{\delta_2}
\Delta_{v'\omega_j, u'\omega_j})^{-a_{ji}}$ of $P'_1$ and $P'_2$
is relatively prime with $\Delta_{s_i \omega_i, \omega_i}$,
it remains to show that
$$P_1 (\Delta_{v'\omega_j, u'\omega_j} \Delta_{v's_j\omega_j, u's_j\omega_j})^{-a_{ji}}
+ P_2 (\Delta_{v'\omega_j, u's_j\omega_j} \Delta_{v's_j \omega_j, u'\omega_j})^{-a_{ji}}$$
is divisible by $\Delta_{s_i \omega_i, \omega_i}$.
Since $P_1 + P_2$ is divisible by $\Delta_{s_i \omega_i,\omega_i}$,
it suffices to show that
$$(\Delta_{v'\omega_j, u'\omega_j} \Delta_{v's_j\omega_j, u's_j\omega_j})^{-a_{ji}}
- (\Delta_{v'\omega_j, u's_j\omega_j} \Delta_{v's_j \omega_j, u'\omega_j})^{-a_{ji}}$$
is divisible by $\Delta_{s_i \omega_i,\omega_i}$.
This in turn follows from the fact that
$$\Delta_{v'\omega_j, u'\omega_j} \Delta_{v's_j\omega_j, u's_j\omega_j}
- \Delta_{v'\omega_j, u's_j\omega_j} \Delta_{v's_j \omega_j, u'\omega_j}$$
is divisible by $\Delta_{s_i \omega_i,\omega_i}$.
But the last expression can be factored according to
(\ref{eq:minors-Dodgson}), and one of the factors is
$\Delta_{s_i \omega_i,\omega_i}^{-a_{ij}}$.
This completes the proofs of Theorem~\ref{th:delta' regular}
and property (3) in Lemma~\ref{lem:main}.

\subsection{Property (4)}
\label{sec:long products}
Let us fix a reduced word $\ii = (i_1, \dots, i_m) \in R(u,v)$,
and an $\ii$-bounded index $n \in [2,m]$.
Let $|i_n| = i \in \Pi$.
Let $\CC^m$ denote the $m$-dimensional vector space with coordinates
$M_1, \dots, M_{n-1}, M'_n, M_{n+1}, \dots, M_m$.
Let $t_1, \dots, t_m$ be rational functions on $\CC^m$ given by
(\ref{eq:negative t-through-x}) and (\ref{eq:positive t-through-x}),
where $M_n$ is determined from (\ref{eq:M'}).
By Theorem~\ref{th:t-through-x}, the map
$$\pi: (M_1, \dots, M_{n-1}, M'_n, M_{n+1}, \dots, M_m) \mapsto
x_\ii (t_1, \dots, t_m)$$
is a birational isomorphism $\CC^m \to L^{u,v}$ inverse to the map
$$x \mapsto (M_1 (x), \dots, M_{n-1} (x), M'_n (x), M_{n+1}(x), \dots,
M_m(x))\ .$$
To prove property (4) in Lemma~\ref{lem:main}, it suffices to show
that $\pi$ restricts to a regular map $\CC_{\neq 0}^m \to L^{u,v}$.

Let us first show that $\pi$ restricts to a regular map $\CC_{\neq 0}^m \to G$.
In view of (\ref{eq:negative t-through-x}) and (\ref{eq:positive t-through-x}),
if $k < n^-$ or $k > n$ then $t_k$ is a Laurent monomial in the
variables $M_l$ with $l \neq n$.
Thus we only need to show that the product
$x_{i_{n^-}}(t_{n^-}) \cdots x_{i_n}(t_n)$ is a regular function
on $\CC_{\neq 0}^m$.
Without loss of generality, we can assume that $n^- = 1$.
For each $k = 2, \dots, n$, we define $p_k \in \CC$ and a rational map
$y_k: \CC^m  \to G$ as follows:
$$p_k = t_1 \doublesubscript{\prod}{1 < l < k}{\veps(i_l) = -1}
t_l^{a_{|i_l|,i}} \ ,$$
while $y_k = x_{i_1} (p_k) x_{i_k}(t_k) x_{i_1} (p_{k+1})^{-1}$ for $k < n$,
and $y_n = x_{i_1} (p_n) x_{i_n}(t_n)$.
Then we have $x_{i_{1}}(t_{1}) \cdots x_{i_n}(t_n) = y_2 \cdots y_n$.
Thus it suffices to show that each $y_k$ is a regular function
on $\CC_{\neq 0}^m$.
As in section \ref{sec:regularity}, we denote
$E_{\pm} = \{l \in [2,n-1] \ : \ \veps (i_l) = \pm 1 \}$.
Let us first prove that $y_n$ is a regular function
on $\CC_{\neq 0}^m$.
We have four cases to consider.

\smallskip

\noindent {\bf Case 1:} $i_1 = i_n = i$.
Then
$$y_n = x_i (p_n) x_i (t_n) = \varphi_i
\mat{1}{p_n}{0}{1} \varphi_i
\mat{1}{t_n}{0}{1} = \varphi_i
\mat{1}{p_n + t_n}{0}{1} \ .$$
It remains to show that $p_n + t_n$ is a regular function on
$\CC_{\neq 0}^m$.
This follows by a direct calculation using
(\ref{eq:negative t-through-x}), (\ref{eq:positive t-through-x})
and (\ref{eq:M'}): we obtain that
$$p_n + t_n = \frac{M'_n}{M_1 M_{n^+}} \cdot
\doublesubscript{\prod}{l > n}{l^- \in E_- \cup \{0\}}
M_l^{C_{ln}} \cdot \doublesubscript{\prod}{l \in E_-}{l^- \in E_+}
M_l^{-C_{ln}}$$
is a Laurent monomial in $M_1, \dots, M_{n-1}, M'_n, M_{n+1}, \dots, M_m$.

\smallskip

\noindent {\bf Case 2:} $i_1 = i,  i_n = -i$.
Then
$$y_n = x_i (p_n) x_{-i} (t_n) = \varphi_i \mat{1}{p_n}{0}{1}
\varphi_i \mat{t_n^{-1}}{0}{1}{t_n} = \varphi_i
\mat{p_n + t_n^{-1}}{p_n t_n}{1}{t_n} \ .$$
It remains to show that each of $t_n$, $p_n t_n$, and $p_n + t_n^{-1}$
is a regular function on $\CC_{\neq 0}^m$.
By a direct calculation, $p_n t_n$ and $p_n + t_n^{-1}$ are
Laurent monomials in $M_1, \dots, M_{n-1}, M'_n, M_{n+1}, \dots,
M_m$, while $t_n = M_n / M_{n^+}$ is the sum of two such
Laurent monomials; in fact, we have
$$p_n + t_n^{-1} = \frac{M'_n}{M_1} \cdot
\doublesubscript{\prod}{l \in E_-}{l^- \in E_+}
M_l^{-C_{ln}} \ .$$

\smallskip

\noindent {\bf Case 3:} $i_1 = -i,  i_n = i$.
Then
$$y_n = x_{-i} (p_n) x_{i} (t_n) = \varphi_i \mat{p_n^{-1}}{0}{1}{p_n}
\varphi_i \mat{1}{t_n}{0}{1}
 = \varphi_i \mat{p_n^{-1}}{p_n^{-1} t_n}{1}{p_n + t_n} \ .$$
It remains to show that each of $p_n^{-1}$, $p_n^{-1} t_n$, and $p_n + t_n$
is a regular function on $\CC_{\neq 0}^m$.
By a direct calculation, $p_n^{-1} t_n$ and $p_n + t_n$ are
Laurent monomials in $M_1, \dots, M_{n-1}, M'_n, M_{n+1}, \dots,
M_m$, while $p_n^{-1}$ is the sum of two such
Laurent monomials; in fact, we have
$$p_n + t_n = \frac{M'_n}{M_{n^+}} \cdot
\doublesubscript{\prod}{l > n}{l^- \in E_-}
M_l^{C_{ln}} \ .$$

\smallskip

\noindent {\bf Case 4:} $i_1 = i_n = -i$.
Then
$$y_n = x_{-i} (p_n) x_{-i} (t_n) = \varphi_i \mat{p_n^{-1}}{0}{1}{p_n}
\varphi_i \mat{t_n^{-1}}{0}{1}{t_n}
 = \varphi_i \mat{p_n^{-1}t_n^{-1}}{0}{p_n + t_n^{-1}}{p_n t_n} \ .$$
It remains to show that each of $(p_n t_n)^{\pm 1}$ and $p_n + t_n^{-1}$
is a regular function on $\CC_{\neq 0}^m$.
By a direct calculation, both $p_n t_n$ and $p_n + t_n^{-1}$ are
Laurent monomials in $M_1, \dots, M_{n-1}, M'_n, M_{n+1}, \dots, M_m$;
in fact, we have
$$p_n + t_n^{-1} = M'_n \cdot
\doublesubscript{\prod}{l \in E_-}{l^- \in E_+ \cup \{0\}}
M_l^{-C_{ln}} \ .$$

\smallskip

Now let $1 < k < n$, and suppose $|i_k| = j \in \Pi$.
To show that $y_k$ is a regular function on $\CC_{\neq 0}^m$,
we have to consider another four cases.

\smallskip

\noindent {\bf Case 1:} $i_1 = i, i_k = j$.
Then $y_k = x_i (p_k) x_j (t_k) x_i (- p_k)$.
Let $p = p_k^{-1}, q = p_k^{- a_{ij}} t_k$.
Clearly, both $p$ and $q$ are regular functions on $\CC_{\neq 0}^m$.
The desired regularity of $y_k$ becomes a consequence of the
following lemma.

\begin{lemma}
\label{lem:rank 2 regularity}
For any two distinct $i,j \in \Pi$, the map
$\CC_{\neq 0} \times \CC \to N$ given by
$(p,q) \mapsto x_i (p^{-1}) x_j (p^{- a_{ij}}q) x_i (- p^{-1})$
extends to a regular map $\CC^2 \to N$.
\end{lemma}

In order not to interrupt the exposition, we will prove this lemma
in the end of this section.

\smallskip

\noindent {\bf Case 2:} $i_1 = i, i_k = -j$.
Then
$$y_k = x_i (p_k) x_{-j} (t_k) x_i (p_k t_k^{a_{ji}})^{-1} = x_{-j}(t_k)
= x_{-j} (M_k / M_{k^+})$$
(see \cite[Proposition~7.2]{BZ99}),
which is a regular function on $\CC_{\neq 0}^m$.

\smallskip

\noindent {\bf Case 3:} $i_1 = -i, i_k = j$.
Then
$$y_k = x_{-i} (p_k) x_{j} (t_k) x_{-i} (p_k)^{-1} = x_{j}(p_k^{a_{ij}} t_k)$$
(see \cite[Proposition~7.2]{BZ99}),
which is a regular function on $\CC_{\neq 0}^m$ since
$p_k^{a_{ij}} t_k$ is a Laurent monomial in
$M_1, \dots, M_{n-1}, M_{n+1}, \dots, M_m$.

\smallskip

\noindent {\bf Case 4:} $i_1 = -i, i_k = -j$.
Then
$$y_k = x_{-i} (p_k) x_{-j} (t_k) x_{-i} (p_k t_k^{a_{ji}})^{-1} \ .$$
Using (\ref{eq:xnegative}) and the commutation relation
\cite[(2.5)]{FZ}, we can rewrite $y_k$ as follows:
$$y_k = y_i (p_k) y_j (p_k^{a_{ij}} t_k) y_i (- p_k)
t_k^{- s_i \alpha_j^\vee} \ .$$
The ``Cartan factor" $t_k^{- s_i \alpha_j^\vee}$ is clearly
a regular function on $\CC_{\neq 0}^m$.
As for the first factor
$y_i (p_k) y_j (p_k^{a_{ij}} t_k) y_i (-p_k)$,
after applying the automorphism  $x \mapsto x^\iota T$ of $G$,
it becomes $x_i (p_k) x_j (p_k^{a_{ij}} t_k) x_i (-p_k)$,
and its regularity follows from Lemma~\ref{lem:rank 2 regularity}
with $p = p_k^{-1}$ and $q = t_k$.

\smallskip

We have proved (modulo Lemma~\ref{lem:rank 2 regularity}) that the map
$$\pi: (M_1, \dots, M_{n-1}, M'_n, M_{n+1}, \dots, M_m) \mapsto
x_\ii (t_1, \dots, t_m)$$
is a regular map $\CC_{\neq 0}^m \to G$.
To complete the proof of property (4), it remains to show that the
image of $\pi$ is contained in $L^{u,v}$.
By Theorem~\ref{th:t-through-x}, this image is contained in the
closure of $L^{u,v}$.
Recall from \cite{FZ} that $L^{u,v}$ is determined
inside its closure by the conditions
$\Delta_{\omega_j, v^{-1} \omega_j} (x) \neq 0$ for all $j$.
The results in \cite{FZ} also imply that, for any $x \in L^{u,v}$,
we have
$\Delta_{\omega_j, v^{-1} \omega_j} (x) =
(\Delta_{\omega_j, u^{-1} \omega_j} (\psi^{u,v}(x)))^{-1}
= M_{k(j)}(x)^{-1}$, where $k(j)$ is the first occurrence of the
index $\pm j$ in $\ii$.
It follows that
$$\Delta_{\omega_j, v^{-1} \omega_j}
(\pi(M_1, \dots, M_{n-1}, M'_n, M_{n+1}, \dots, M_m)) =
M_{k(j)}^{-1} \neq 0$$
on $\CC_{\neq 0}^m$, and we are done.

\smallskip

\noindent {\sl Proof of Lemma~\ref{lem:rank 2 regularity}.}
The proof below was suggested by N.~Reshetikhin; it is simpler
than the author's original proof.
Let $\{e_i \ : i \in \Pi\}$ be the standard generators of the Lie
algebra $\nn = {\rm Lie} (N)$; thus, we have
$x_i (t) = \exp (t e_i)$ for every $i \in \Pi$ and $t \in \CC$.
Let ${\rm Ad}  : G \to {\rm Aut} (\gg)$ be the adjoint
representation of $G$, and ${\rm ad}  : \gg \to {\rm End} (\gg)$
be the differential of ${\rm Ad}$; recall that these
representations satisfy
$\exp \ ({\rm Ad}(x) e) = x \exp \ (e) x^{-1}$ and
$\exp \ ({\rm ad}(e)) = {\rm Ad} (\exp \ (e))$ for $x \in G$ and
$e \in \gg$.
It follows that
$$x_i (p^{-1}) x_j (p^{- a_{ij}}q) x_i (- p^{-1}) =
\exp \ ({\rm Ad}(x_i (p^{-1})) \cdot p^{- a_{ij}}q e_j)$$
$$= \exp \ ( \exp \ ({\rm ad} (p^{-1} e_i)) \cdot p^{- a_{ij}}q e_j)
= \exp \ \left(q \sum_{n \geq 0} \frac{p^{- a_{ij}-n}}{n!}
{\rm ad} (e_i)^n (e_j)\right)$$
$$ = \exp \ \left(q \sum_{n = 0}^{- a_{ij}} \frac{p^{- a_{ij}-n}}{n!}
{\rm ad} (e_i)^n (e_j)\right) \ ,$$
which is obviously regular in $p$ and $q$ (the last equality
follows from Serre's relation ${\rm ad} (e_i)^{1 - a_{ij}} (e_j) = 0$).
\endproof

Property (4) in Lemma~\ref{lem:main}, and Theorems~\ref{th:components explicit}
and \ref{th:components} are finally proved.

\section{Some examples and applications}
\label{sec:examples}

\subsection{Cones of regular monomials}
\label{sec:regular t-monomials}
Let us again fix a pair $(u,v) \in W \times W$, and a reduced word
$\ii = (i_1, \ldots, i_m) \in R(u,v)$.
In view of Theorem~\ref{th:t-through-x}, a generic element $x \in L^{u,v}$
has the form $x= x_\ii (t_1, \ldots, t_m)$,
so the factorization parameters $t_k$ are well-defined rational functions on $L^{u,v}$
given by (\ref{eq:negative t-through-x}) and (\ref{eq:positive t-through-x}).
Combining these formulas with Lemma~\ref{lem:regular monomials} yields the following corollary.

\begin{proposition}
\label{pr:regular t-monomials}
A Laurent monomial $t_1^{a_1} \cdots t_m^{a_m}$
is a regular function on $L^{u,v}$ if and only if
\begin{equation}
\label{eq:t-regular forms}
- \veps (i_n) a_n - a_{n^-} + \doublesubscript{\sum}
{n^- < k <  n}{\veps (i_k) = 1} C_{nk} a_k \geq 0
\end{equation}
for any $\ii$-bounded $n \in [1,m]$.
\end{proposition}

Two special cases are worth mentioning.
If $v = e$ then $\veps (i_k) = -1$ for all $k$, and the
inequalities (\ref{eq:t-regular forms}) take the form
$a_n \geq a_{n^-}$.
If $u = e$ then $\veps (i_k) = 1$ for all $k$, and the
inequalities (\ref{eq:t-regular forms}) take the form
$$- a_n - a_{n^-} + \sum_{n^- < k <  n} C_{nk} a_k \geq 0
\ ;$$
the cone defined by these inequalities appeared in a different
context in \cite{L}, and also in \cite{NZ}.

\subsection{Intersections of two open opposite Schubert cells}
\label{sec:L e,wnot}
Let us illustrate Theorem~\ref{th:components} by the case
when $u = e$ and $w = w_0$, the longest element in $W$.
In this case, $L^{u,v}$ is biregularly
isomorphic to the intersection of two open opposite Schubert cells
$C_{w_0} \cap w_0 C_{w_0}$,
where $C_{w_0} = (B w_0 B)/B$ is the open Schubert cell
in the flag variety $G/B$.
These opposite cells appeared in the literature in various
contexts, and were studied (in various degrees of generality)
in \cite{BFZ,BZ97,RI,RI2,SSV1,SSV2}.
Let $C$ denote the number of connected components of
$L^{e,w_0}(\RR)$; to emphasize the dependency on $G$, we shall
write $C = C(X_r)$, where $X_r = A_r, B_r, \dots, G_2$ is the type
of $G$ in the Cartan-Killing classification.

The numbers $C(A_r)$ were determined in \cite{SSV1,SSV2}: it turns
out that $C(A_1) = 2, C(A_2) = 6, C(A_3) = 20, C(A_4)= 52$, and
$C(A_r) = 3 \cdot 2^{r}$ for $r \geq 5$.
Theorem~\ref{th:components} allows us to extend this result to all
other simply-laced types.

\begin{proposition}
\label{pr:DE count}
If $X_r$ is one of the types $A_r \  (r \geq 5), D_r \  (r \geq 4), E_6, E_7$, or
$E_8$ then $C(X_r) = 3 \cdot 2^{r}$.
\end{proposition}

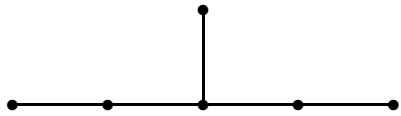
\begin{figure}[ht]
\setlength{\unitlength}{1.8pt}

\begin{center}
\begin{picture}(80,20)(0,0)
\thicklines

\put(0,0){\line(1,0){80}}
\put(40,0){\line(0,1){20}}

  \put(0,0){\circle*{2.0}}
  \put(20,0){\circle*{2.0}}
  \put(40,0){\circle*{2.0}}
  \put(60,0){\circle*{2.0}}
  \put(80,0){\circle*{2.0}}
  \put(40,20){\circle*{2.0}}

\end{picture}
\end{center}

\caption{The Dynkin graph $E_6$.}
\label{fig:E6}
\end{figure}

\proof
Following \cite[Definition~3.10]{SSVZ}, we say that
a graph is $E_6$-compatible if it is connected,
and it contains an induced subgraph with $6$ vertices isomorphic to the Dynkin
graph $E_6$ (see Fig.~\ref{fig:E6}).

Combining Theorem~\ref{th:components} with \cite[Corollary~3.12]{SSVZ},
we obtain the following sufficient condition for the equality
$C(X_r) = 3 \cdot 2^{r}$: it holds provided $G$ is simply-laced,
and there exists $\ii \in R(w_0)$ such that
the induced subgraph of $\Sigma(\ii)$ (see Definition~\ref{def:edges})
on the set of all $\ii$-bounded vertices is $E_6$-compatible.
In \cite{SSV2} this condition was checked for the type $A_5$.
Therefore, it also holds for any simply-laced Dynkin graph that
contains an induced subgraph of type $A_5$, that is, for
$A_r \ (r \geq 5), D_r \  (r \geq 6), E_6, E_7$, and $E_8$.
It remains to check this condition for the type $D_4$ (the
statement for $D_5$ then follows).
Let $\Pi = \{1,2,3,4\}$ with the branching vertex $3$.
Take the reduced word $\ii = (1,2,3,1,2,3,4,3,1,2,3,4) \in R(w_0)$.
By inspection, the induced subgraph of $\Sigma (\ii)$ with
$\ii$-bounded vertices $4,5,9,10,11$, and $12$ is isomorphic to the Dynkin
graph $E_6$, and we are done.
\endproof

The numbers $C(B_2)$ and $C(G_2)$ were determined in \cite{RI2}:
it turns out that $C(B_2)= 8$ and $C(G_2) = 11$.
Theorem~\ref{th:components} gives a simpler way to prove these
answers.
In the case of $B_2$, take $\ii = (j,i,j,i)$ with $a_{ij} = -2$
and $a_{ji} = -1$.
Then $\Gamma_\ii (\FF_2)$ is the group of transformations of
$\FF_2^4$ generated by $\tau_3 : \xi_3 \to \xi_3 + \xi_1$ and
$\tau_4 : \xi_4 \to \xi_4 + \xi_3 + \xi_2$.
It is easy to see that the action of $\Gamma_\ii (\FF_2)$
in $\FF_2^4$ has $8$ orbits: four fixed points $0000, 0001, 0110$, and $0111$,
two $2$-element orbits $0010 \overset {\tau_4}{\longleftrightarrow} 0011$
and $0100 \overset {\tau_4}{\longleftrightarrow} 0101$,
and two $4$-element orbits $1000 \overset {\tau_3}{\longleftrightarrow} 1010
\overset {\tau_4}{\longleftrightarrow} 1011 \overset {\tau_3}{\longleftrightarrow} 1001$
and $1110 \overset {\tau_3}{\longleftrightarrow} 1100
\overset {\tau_4}{\longleftrightarrow} 1101 \overset {\tau_3}{\longleftrightarrow}
1111$.

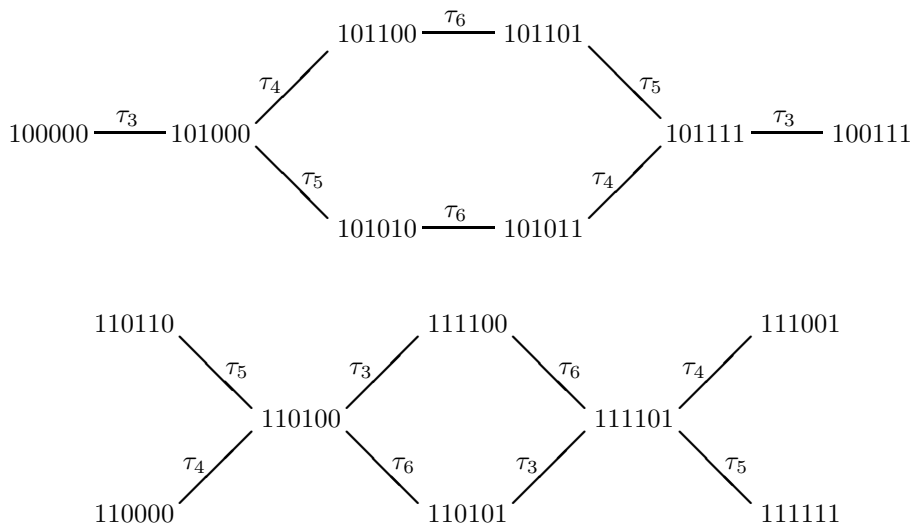
\begin{figure}[ht]
\setlength{\unitlength}{1.8pt}

\begin{center}
\begin{picture}(160,105)(0,0)
\thicklines

\put(0,80){\line(1,0){15}}
\put(34,82){\line(1,1){15}}
\put(34,77){\line(1,-1){15}}
\put(69,101){\line(1,0){15}}
\put(69,60){\line(1,0){15}}
\put(104,98){\line(1,-1){15}}
\put(104,62){\line(1,1){15}}
\put(138,80){\line(1,0){15}}

\put(-18,78){${100000}$}
\put(16,78){${101000}$}
\put(51,99){${101100}$}
\put(51,58){${101010}$}
\put(86,99){${101101}$}
\put(86,58){${101011}$}
\put(120,78){${101111}$}
\put(155,78){${100111}$}

\put(7,83){\makebox(0,0){$\tau_3$}}
\put(145,83){\makebox(0,0){$\tau_3$}}
\put(76,104){\makebox(0,0){$\tau_6$}}
\put(76,63){\makebox(0,0){$\tau_6$}}
\put(37,90){\makebox(0,0){$\tau_4$}}
\put(107,70){\makebox(0,0){$\tau_4$}}
\put(117,90){\makebox(0,0){$\tau_5$}}
\put(46,70){\makebox(0,0){$\tau_5$}}

\put(18,37){\line(1,-1){15}}
\put(18,2){\line(1,1){15}}
\put(53,22){\line(1,1){15}}
\put(53,17){\line(1,-1){15}}
\put(88,37){\line(1,-1){15}}
\put(88,2){\line(1,1){15}}
\put(123,22){\line(1,1){15}}
\put(123,17){\line(1,-1){15}}

\put(0,38){${110110}$}
\put(0,-2){${110000}$}
\put(35,18){${110100}$}
\put(70,38){${111100}$}
\put(70,-2){${110101}$}
\put(105,18){${111101}$}
\put(140,38){${111001}$}
\put(140,-2){${111111}$}

\put(30,30){\makebox(0,0){$\tau_5$}}
\put(21,10){\makebox(0,0){$\tau_4$}}
\put(56,30){\makebox(0,0){$\tau_3$}}
\put(65,10){\makebox(0,0){$\tau_6$}}
\put(91,10){\makebox(0,0){$\tau_3$}}
\put(100,30){\makebox(0,0){$\tau_6$}}
\put(126,30){\makebox(0,0){$\tau_4$}}
\put(135,10){\makebox(0,0){$\tau_5$}}

\end{picture}
\end{center}

\caption{The $8$-element orbits for $G_2$.}
\label{fig:8-orbits}
\end{figure}

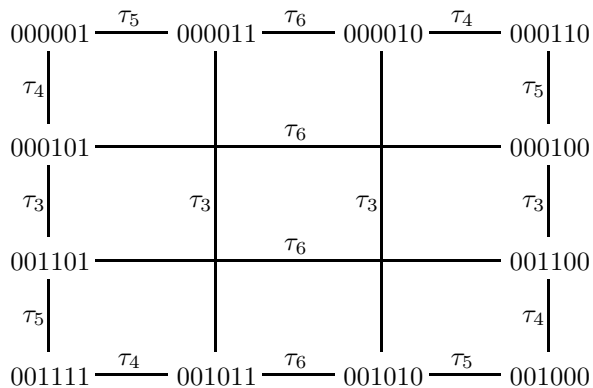
\begin{figure}[ht]
\setlength{\unitlength}{1.8pt}

\begin{center}
\begin{picture}(120,80)(0,0)
\thicklines

\put(0,-2){${001111}$}
\put(35,-2){${001011}$}
\put(70,-2){${001010}$}
\put(105,-2){${001000}$}
\put(0,22){${001101}$}
\put(105,22){${001100}$}
\put(0,46){${000101}$}
\put(105,46){${000100}$}
\put(0,70){${000001}$}
\put(105,70){${000110}$}
\put(35,70){${000011}$}
\put(70,70){${000010}$}

\put(18,0){\line(1,0){15}}
\put(53,0){\line(1,0){15}}
\put(88,0){\line(1,0){15}}
\put(8,5){\line(0,1){15}}
\put(113,5){\line(0,1){15}}
\put(8,29){\line(0,1){15}}
\put(113,29){\line(0,1){15}}
\put(8,53){\line(0,1){15}}
\put(113,53){\line(0,1){15}}
\put(18,72){\line(1,0){15}}
\put(53,72){\line(1,0){15}}
\put(88,72){\line(1,0){15}}
\put(18,24){\line(1,0){85}}
\put(18,48){\line(1,0){85}}
\put(43,5){\line(0,1){63}}
\put(78,5){\line(0,1){63}}

\put(25,3){\makebox(0,0){$\tau_4$}}
\put(60,3){\makebox(0,0){$\tau_6$}}
\put(95,3){\makebox(0,0){$\tau_5$}}
\put(5,12){\makebox(0,0){$\tau_5$}}
\put(110,12){\makebox(0,0){$\tau_4$}}
\put(5,36){\makebox(0,0){$\tau_3$}}
\put(110,36){\makebox(0,0){$\tau_3$}}
\put(5,60){\makebox(0,0){$\tau_4$}}
\put(110,60){\makebox(0,0){$\tau_5$}}
\put(25,75){\makebox(0,0){$\tau_5$}}
\put(60,75){\makebox(0,0){$\tau_6$}}
\put(95,75){\makebox(0,0){$\tau_4$}}
\put(40,36){\makebox(0,0){$\tau_3$}}
\put(75,36){\makebox(0,0){$\tau_3$}}
\put(60,27){\makebox(0,0){$\tau_6$}}
\put(60,51){\makebox(0,0){$\tau_6$}}

\end{picture}
\end{center}

\caption{The $12$-element orbit for $G_2$.}
\label{fig:12-orbit}
\end{figure}

The case of $G_2$ is treated in a similar fashion.
Take $\ii = (j,i,j,i,j,i)$ with $a_{ij} = -3$ and $a_{ji} = -1$.
Then $\Gamma_\ii (\FF_2)$ is the group of transformations of
$\FF_2^6$ generated by the transvections $\tau_n \ (3 \leq n \leq 6)$
acting by $\tau_n: \xi_n \to \sum_{|k-n| \leq 2} \xi_k$.
It is easy to see that the action of $\Gamma_\ii (\FF_2)$
in $\FF_2^6$ has $11$ orbits: four fixed points $000000, 001001, 001110$, and
$000111$; six $8$-element orbits, and one $12$-element orbit.
The $8$-element orbits are depicted in Fig.~\ref{fig:8-orbits}
(one has to take the first depicted orbit together with its
translates by the $3$ non-zero fixed vectors; and the second
depicted orbit together with its translate by the vector
$001110$); the $12$-element orbit is depicted in
Fig.~\ref{fig:12-orbit}.

\begin{remark}
{\rm Computing the numbers $C(B_r)$ and $C(C_r)$ for $r \geq 3$ seems
to be a challenging problem.
Since the transvections $\tau_n$ are no longer symplectic in this
case, one cannot use \cite[Corollary~3.12]{SSVZ}
(at least, not in a straightforward way).}
\end{remark}

\subsection{Dual canonical basis for the type $B_2$}
\label{sec:B2 basis}
In conclusion, we briefly discuss a potential application of the above
results.
Let $G/N$ be the \emph{base affine space} for $G$.
It is well-known that the ring of regular functions $\CC[G/N]$ (that is,
regular functions on $G$ invariant under right translations by
elements of $N$) is the multiplicity-free sum of all irreducible
finite-dimensional representations of $G$.
Let $B$ denote the \emph{dual canonical basis} in $\CC[G/N]$
(more precisely, $B$ is the ``classical limit" of the dual canonical basis
in the $q$-deformed ring $\CC_q [G/N]$).
Despite much progress in studying properties of the canonical bases, an explicit
construction of $B$ still remains to be found.
It is known that $B$ contains all ``Pl\"ucker coordinates"
$P_\gamma = \Delta_{\gamma, \omega_i}$
for $i \in \Pi$ and $\gamma \in W \omega_i$.
We suspect that $B$ also contains all functions
$\Delta'_\ii$ in Theorem~\ref{th:delta' regular}
corresponding to reduced words $\ii$ consisting of elements of $\Pi$.
Thus, these functions together with the Pl\"ucker coordinates
$P_\gamma$ are among the building blocks for $B$.

As an illustration, consider the case when $G$ is of type
$B_2$, i.e., $\Pi = \{i,j\}$ with $a_{ij} = -2$ and $a_{ji} = -1$.
The basis $B$ in this case was found in \cite{RZ}
(even before the ``official" discovery of canonical bases).
Translating the results in \cite{RZ} into our present notation,
we obtain the following.

There are $8$ Pl\"ucker coordinates: $P_{\omega_i}, P_{\omega_j},
P_{s_i \omega_i}, P_{s_j \omega_j}, P_{s_j s_i \omega_i},
P_{s_i s_j \omega_j}, P_{w_0 \omega_i}$, and $P_{w_0 \omega_j}$.
Let us also denote $Q_{\omega_j} = \Delta'_{(i,j,i)}$ and
$Q_{2 \omega_i} = \Delta'_{(j,i,j)}$; thus, these functions are
defined from the equations
\begin{equation}
\label{eq:Q omegaj}
Q_{\omega_j} P_{s_i \omega_i} = P_{s_i s_j \omega_j} P_{\omega_i}
+ P_{\omega_j} P_{w_0 \omega_i}
\end{equation}
and
\begin{equation}
\label{eq:Q omegai}
Q_{2 \omega_i} P_{s_j \omega_j} = P_{s_j s_i \omega_i}^2  P_{\omega_j}
+ P_{\omega_i}^2 P_{w_0 \omega_j} \ .
\end{equation}
The main result of \cite{RZ} can be now summarized as follows.

\begin{figure}[ht]
\setlength{\unitlength}{1.8pt}

\begin{center}
\begin{picture}(90,40)(0,0)
\thicklines

\put(0,22){$Q_{2 \omega_i}$}
\put(21,44){$P_{s_i \omega_i}$}
\put(21,-2){$P_{s_j s_i \omega_i}$}
\put(55,44){$P_{s_i s_j \omega_j}$}
\put(59,-2){$P_{s_j \omega_j}$}
\put(80,22){$Q_{\omega_j}$}

\put(8,26){\line(1,1){15}}
\put(8,19){\line(1,-1){15}}
\put(33,46){\line(1,0){20}}
\put(37,0){\line(1,0){20}}
\put(65,4){\line(1,1){15}}
\put(65,41){\line(1,-1){15}}

\end{picture}
\end{center}

\caption{The ``magical hexagon" for $B_2$.}
\label{fig:hexagon for B2}
\end{figure}
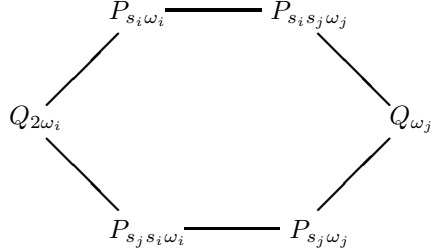

\begin{proposition}
\label{pr:B2 basis}
The dual canonical basis $B$ of $\CC[G/N]$ consists of all
monomials in $10$ variables $P_{\omega_i}, \dots, P_{w_0 \omega_j},
Q_{\omega_j}, Q_{2 \omega_i}$ with the following property: if this
monomial contains variables in two vertices of the ``magical
hexagon" in Fig.~\ref{fig:hexagon for B2} then these two vertices
are adjacent.
\end{proposition}

We see that $B$ is the union (not disjoint!) of six families
of elements corresponding to the edges of the hexagon in
Fig.~\ref{fig:hexagon for B2}: each family consists of all
monomials in six variables $P_{\omega_i}, P_{\omega_j},
P_{w_0 \omega_i}, P_{w_0 \omega_j}, P, Q$, where
$P$ and $Q$ lie in two adjacent vertices of the hexagon.

Note that the equations (\ref{eq:Q omegaj}) and (\ref{eq:Q omegai})
can be now interpreted as expansions in the basis $B$ of two
``forbidden" monomials corresponding to diagonals of the hexagon.
There are 7 more such identities corresponding to the remaining 7
diagonals:
$$Q_{\omega_j} P_{s_j s_i \omega_i} = P_{w_0 \omega_j} P_{\omega_i}
+ P_{s_j \omega_j} P_{w_0 \omega_i} \ ;$$

$$Q_{2 \omega_i} P_{s_i s_j \omega_j} = P_{w_0 \omega_i}^2  P_{\omega_j}
+ P_{s_i \omega_i}^2 P_{w_0 \omega_j} \ ;$$

$$P_{s_i \omega_i} P_{s_j s_i \omega_i} = P_{\omega_i} P_{w_0 \omega_i}
+ Q_{2 \omega_i} \ ;$$

$$P_{s_j \omega_j} P_{s_i s_j \omega_j} = P_{\omega_j} P_{w_0 \omega_j}
+ Q_{\omega_j}^2 \ ;$$

$$P_{s_i \omega_i} P_{s_j \omega_j} = P_{s_j s_i \omega_i} P_{\omega_j}
+ P_{\omega_i} Q_{\omega_j} \ ;$$

$$P_{s_j s_i \omega_i} P_{s_i s_j \omega_j} = P_{s_i \omega_i} P_{w_0 \omega_j}
+ P_{w_0 \omega_i} Q_{\omega_j} \ ;$$

$$Q_{2 \omega_i} Q_{\omega_j} = P_{s_j s_i \omega_i} P_{w_0 \omega_i}
 P_{\omega_j} +  P_{s_i \omega_i} P_{\omega_i}
 P_{w_0 \omega_j} \ .$$

\section*{Acknowledgements}
I thank Arkady Berenstein, Nicolai Reshetikhin, Boris and Misha Shapiro, and  Alek
Vainshtein for helpful discussions.

\end{document}